\documentclass[a4paper, 11pt]{article}
\usepackage[utf8]{inputenc}
\usepackage[english]{babel}
\usepackage{expl3}
\usepackage{array}
\usepackage{url}
\usepackage{amsmath,amsthm,amsfonts,amssymb,amscd}
\usepackage{fullpage}
\usepackage{algorithm}
\usepackage{algpseudocode}
\usepackage{tabularx}
\usepackage{booktabs}
\usepackage{lipsum}
\usepackage{bm}
\usepackage[colorlinks,urlcolor=cobalt,citecolor=cobalt,linkcolor=cobalt,pdftex, pdfauthor = Lukas Exl]{hyperref}
\usepackage{siunitx}
\usepackage{graphicx}
\usepackage{xcolor}
\usepackage{caption}
\usepackage{subcaption}
\usepackage[section]{placeins}
\usepackage{listings}
\usepackage{authblk}
\usepackage{stmaryrd} 

\usepackage{url}

\definecolor{cobalt}{rgb}{0.0, 0.28, 0.67}

\theoremstyle{definition}

\DeclareUnicodeCharacter{2113}{\ensuremath{\ell}}

\renewcommand{\vec}[1]{\boldsymbol{#1}}

\makeatletter

\newcounter{phase}[algorithm]
\newlength{\phaserulewidth}
\newcommand{\setphaserulewidth}{\setlength{\phaserulewidth}}

\makeatother
\setphaserulewidth{.7pt}

\DeclareSymbolFont{largesymbolsA}{U}{txexa}{m}{n}
\DeclareMathSymbol{\varprod}{\mathop}{largesymbolsA}{16}

\date{}

\title{Physics-informed low-rank neural operators with application to parametric elliptic PDEs}
\author[a,b]{Sebastian Schaffer \thanks{\texttt{sebastian.schaffer@univie.ac.at}}}
\author[a,b,c]{Lukas Exl} 

\affil[a]{Wolfgang Pauli Institute, Vienna, Austria}
\affil[b]{Research Platform MMM Mathematics-Magnetism-Materials, University of Vienna, Vienna, Austria}
\affil[c]{Faculty of Mathematics, University of Vienna, Vienna} 

\begin{document}
\maketitle
\noindent\textbf{Abstract.}
\noindent We present the Physics-Informed Low-Rank Neural Operator (PILNO), a neural operator
framework for efficiently approximating solution operators of partial differential equations (PDEs) on point cloud data. PILNO combines low-rank kernel approximations with an encoder--decoder architecture, enabling fast, continuous one-shot predictions while remaining independent of specific discretizations. The model is trained using a physics-informed penalty framework, ensuring that PDE constraints and boundary conditions are satisfied in both supervised and unsupervised settings. We demonstrate its effectiveness on diverse problems, including function fitting, the Poisson equation, the screened Poisson equation with variable coefficients, and parameterized Darcy flow. The low-rank structure provides computational efficiency in high-dimensional parameter spaces, establishing PILNO as a scalable and flexible surrogate modeling tool for PDEs.


\section{Introduction}

Parametric partial differential equations (PDEs) arise throughout science and engineering, for example when material properties, source terms, or boundary conditions depend on parameters. Classical solvers such as finite element or finite difference methods provide highly accurate solutions for a single parameter instance at near-linear cost. However, when solutions are required for many parameter configurations—as in optimization, uncertainty
quantification, or inverse problems—the repeated assembly and solution of large linear systems becomes computationally prohibitive. This challenge worsens in high-dimensional parameter spaces, where the number of samples needed to explore the solution manifold grows exponentially, a phenomenon commonly referred to as the \emph{curse of dimensionality}.

Machine learning methods offer an attractive alternative by enabling surrogate models that approximate the PDE solution (operator) directly. Physics-Informed Neural Networks (PINNs) \cite{raissi2019physics} have gained popularity for approximating PDE solutions, but they suffer from slow and sometimes unstable training, limiting their efficiency in practice. Recent extensions such as \emph{conditional PINNs} \cite{kovacs2022conditional} address this
issue by incorporating parametric dependence directly into the network, allowing families of parametric PDEs to be represented by a single model. While promising, these approaches still face difficulties in scalability and training cost.

A more recent paradigm is operator learning, where the goal is to approximate the solution operator of a PDE directly \cite{lu2021learning, li2024physics}. Once trained, a neural operator provides one-shot predictions for unseen inputs, amortizing the cost of training across potentially large numbers of PDE instances. This perspective connects naturally to convolution-like integral formulations: many PDE solutions can be expressed as
\begin{equation}\label{eq:conv}
(k * f)( x) = \int_\Omega k(x, y)\, f( y)\,\mathrm{d} y,
\end{equation}
where $k$ is a kernel, such as a Green’s function. Approximating such convolutional operators efficiently is central to the success of neural operator architectures.

Neural operators aim to approximate a (possibly nonlinear) operator $\mathcal{T}(f)$ by lifting the input function $f$ into a latent representation and applying iterative convolution-like updates. Specifically, we set $v_0 = f$ (or an embedding of $f$), and update the latent features according to
\begin{equation}\label{eq:iter-scheme}
    v_{t+1}(x) = \sigma\left(W_t v_t(x) +
        \int_\Omega k_t(x, y)\, v_t(y)\,\mathrm{d}y \right),
    \quad t=0,\dots,T,
\end{equation}
where $k_t$ is a trainable kernel, $W_t$ a weight matrix, and $\sigma$ a nonlinear activation function. After $T$ iterations, the final latent representation $v_T$ is projected to the solution space via a neural network $\mathcal{N}_d$, so that $\mathcal{T}(f)(x) \approx \mathcal{N}_d(v_T(x))$.
To reduce the cost of the integral operator, many modern architectures impose a \emph{low-rank} kernel parameterization $k_t(x,y) = \psi_t(x)^\top \phi_t(y)$, which enables efficient evaluation through matrix products on point clouds.

Several architectures instantiate~\eqref{eq:iter-scheme} differently. Convolutional neural networks and DeepONets \cite{lu2021learning} are effective on structured domains, while physics-informed DeepONets can further reduce data requirements \cite{wang2021learning}. Graph-based operators \cite{li2020neural} extend neural operators to unstructured meshes at significant computational cost. Fourier Neural Operators (FNOs) \cite{li2020fourier} achieve
impressive accuracy and efficiency on equispaced grids but require adaptations for general geometries \cite{li2023fourier}. Multi-grid tensorized FNOs \cite{kossaifi2023multi} have been proposed to improve scalability. More recently, low-rank neural operators (LNOs) \cite{kovachki2023neural, zeng2025point} leverage the factorization above to approximate convolutional operators with reduced complexity, which is particularly well suited to
unstructured data and parametric settings.

We propose the \emph{Physics-Informed Low-Rank Neural Operator} (PILNO), which combines low-rank kernel approximations with an encoder–decoder architecture. Unlike grid-based methods such as FNOs, PILNO operates directly on point cloud inputs without requiring structured discretizations. Compared to PINNs and conditional PINNs, it avoids retraining for each parameter configuration by learning the operator itself. In contrast to graph-based methods, PILNO achieves computational efficiency through the low-rank kernel structure. These features position PILNO as a scalable surrogate modeling framework for parametric PDEs.

Once trained, the model delivers continuous one-shot predictions for a given problem. The encoder $E$ maps the functional input into a latent space producing $Z$, which the decoder $D$ then maps to arbitrary points in the domain:
\begin{equation}\label{eq:encdec}
    Z = E(\vec X, f), \qquad u(y) = D(Z, y),
\end{equation}
where $\vec X$ denotes the point cloud, $f$ the input function, and $u(y)$ the predicted solution at point $y$. The computationally heavy iterative updates occur only during encoding, allowing the decoder to remain lightweight and fast. We train PILNO with a physics-informed penalty framework in both supervised and unsupervised settings and demonstrate its effectiveness on parametric elliptic PDEs, highlighting efficiency in high-dimensional parameter spaces where traditional solvers become intractable.

\section{LNO encoder-decoder architecture}
Our goal is to learn a nonlinear approximation of a solution operator for a class of PDEs, without imposing discretization constraints typical of traditional methods. 
We propose an encoder that aggregates information from the input function on a point cloud into a latent representation. A low-rank kernel decomposition makes the integral operator efficient, and the decoder projects this representation to arbitrary query points.

We operate directly on point cloud data in $d$ dimensions: $\vec X = \{x_1,\dots,x_N\} \subset \Omega\subseteq\mathbb{R}^d$, referred to as sensor points, represented as an $N\times d$ array. The model should maintain strong inductive bias: increasing the number of sensor points should not degrade performance as long as the distribution of $\vec X$ is largely preserved. For simplicity, we denote row-wise evaluation as $g(\vec X)$ for $g:\mathbb{R}^d\rightarrow\mathbb{R}^{d'}$, and omit network parameters.

First, the input $f$ is mapped to an $S$-dimensional latent representation $v_0\in \mathbb{R}^S$ via an MLP $\mathcal{N}_0$ applied to all sensor points, $v_0(x_i) = \mathcal{N}_0(f(x_i))$, resulting in $v_0(\vec X)=\mathcal{N}_0(f(\vec X))$. We then iteratively update the latent representation using a more expressive network $\mathcal{N}_t$ per layer:
\begin{equation}\label{eq:iter-scheme2}
    v_{t}(x) = \mathcal{N}_t\Big(LN_t\big(v_{t-1} + \int_\Omega k_t(x, y) v_{t-1}(y)\mathrm{d}y\big)\Big),
\end{equation}
where $LN_t$ is layer normalization. Except for the input and output layers, mappings $\mathbb{R}^S\rightarrow\mathbb{R}^S$ preserve feature dimensions. The low-rank kernel is $k_t(x, y) = \psi_t(x)^T \phi_t(y)$, with $\psi_t, \phi_t: \mathbb{R}^d\to\mathbb{R}^R$. On the point cloud, we approximate the integral using Monte Carlo:
\begin{equation}
    \int_\Omega k_t(x, y) v_{t-1}(y)\mathrm{d}y \approx \frac{|\Omega|}{N} \sum_{i=1}^N \psi(x)^T \phi_t(x_i)v_{t-1}(x_i) = \frac{|\Omega|}{N}\psi_t(x)^T \phi_t(\mathbf{X})^T v_{t-1}(\mathbf{X}).
\end{equation}

The encoder computes:
\begin{equation}\label{eq:encoder}
    \begin{aligned}
        v_0(\mathbf{X}) &= \mathcal{N}_0(f(\mathbf{X})) \\
        h_t(\mathbf{X}) &= v_{t-1}(\mathbf{X}) + \frac{|\Omega|}{N} \psi_t(\mathbf{X}) \phi_t(\mathbf{X})^T v_{t-1}(\mathbf{X}) \\
        v_t(\mathbf{X}) &= \mathcal{N}_t\big(\text{LN}_t(h_t(\mathbf{X}))\big), \quad t=1,\dots,T-1 \\
        z_T(\mathbf{X}) &= \phi_T(\mathbf{X})^T v_{T-1}(\mathbf{X}),
    \end{aligned}
\end{equation}

The decoder computes:
\begin{equation}\label{eq:decoder}
    \begin{aligned}
        v_T(\mathbf{Y}) &= \psi_T(\mathbf{Y}) z_T(\mathbf{X})\\
        u_{\text{LNO}}(\mathbf{Y}) &= \mathcal{N}_T(v_T(\mathbf{Y})) = \mathcal{M}(\vec X, f, \Theta)(\vec Y),
    \end{aligned}
\end{equation}

The encoder transforms $f$ on sensor points into $z_T\in \mathbb{R}^{R\times S}$; the decoder maps this representation to the target points to predict the solution. Kernel networks $\psi_t, \phi_t$ and the final decoder use linear outputs. Monte Carlo can be replaced with other quadrature rules without breaking the low-rank structure. Encoding has complexity $\mathcal{O}(NRST)$ and decoding $\mathcal{O}(MRS)$, scaling linearly with the number of sensor or target points, making the model efficient.

This model leaves substantial room for extensions. The initial mapping could incorporate additional features of $f$, such as its Jacobian $J_f(x_i)$. Notably, the kernel factors $\psi_t$ and $\phi_t$ are independent of $f$, enabling precomputation for a given problem. Kernels may also embed extra information, including geometric features via (approximate) signed distance functions \cite{sukumar2022exact, schaffer2024constraint} or PDE-specific parameters.

For symmetric kernels, one can set $\vec \psi_t = \vec \phi_t$. Importantly, the output $u(x)$ depends solely on the decoder, which simplifies physics-informed training: derivatives with respect to target points are required only for the decoder. Low-rank convolution evaluation scales as $\mathcal{O}(NRS)$ and can be computed efficiently with matrix operations for moderate $R$ and $S$. Encoding with $T$ layers costs $\mathcal{O}(NRST)$, while decoding $M$ target points costs $\mathcal{O}(MRS)$.

Overall, the computational complexity is essentially linear in the number of sensor or target points. Comparing with conventional mesh-based operators—also linear—reveals the strong impact of problem dimensionality. While discretization in high-dimensional parameter spaces is prohibitive, learning an operator allows us to overcome the curse of dimensionality. We present an example in the Results section. Direct comparisons are challenging due to factors such as assembly cost, potential geometry-aware LNOs, and uncertain training strategies. Nonetheless, LNOs likely offer superior speed, although with lower accuracy than conventional mesh-based solvers at present. If high accuracy is required, the LNO model could also be used to compute a “cheap” initial guess for a more accurate grid based method.


\section{Training}
When a large amount of labeled data is available, for instance from simulations, supervised learning with stochastic gradient descent can be applied. Our primary interest lies in unsupervised training tailored to a specific task. The encoder \eqref{eq:encoder} requires an input function $f$ sampled from a functional space $\mathcal{F}$. The model will perform well only if the input function lies approximately within this space. To achieve this, we define a basis for $\mathcal{F}$ and sample its coefficients.

\subsection{B-spline sampling}
We adopt a B-spline tensor product extreme learning machine as in \cite{exl2025higher}. This basis has several advantages: (i) B-splines provide excellent approximation of smooth functions, (ii) coefficient values directly limit corresponding function values, and (iii) the Cox–de Boor recursion ensures efficient and stable evaluation that scales with the B-spline degree due to sparsity.
\begin{figure}
    \centering
    \includegraphics[width=0.7\linewidth]{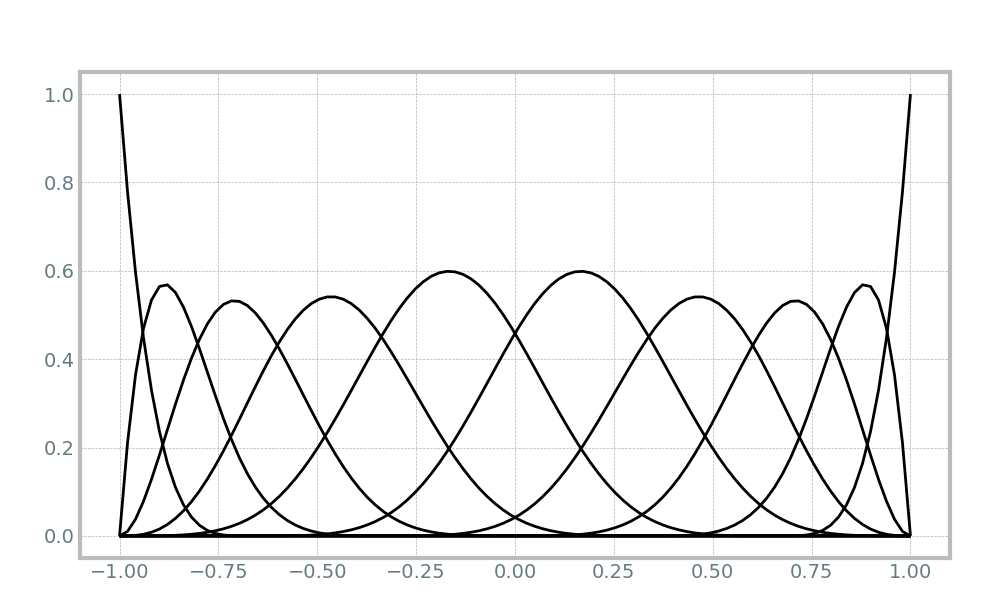}
    \caption{1d B-spline basis of order $K=3$ with $r=7$ equidistant knots resulting in $r + K = 10$ basis functions.}
    \label{fig:B-spline}
\end{figure}

We consider tensor product B-splines $\rho(x) = \rho_1(x^{(1)})\dotsc\rho_d(x^{(d)})$ of order $K$ (degree $K-1$), defined on $r_i$ equidistant knots per direction with $K+1$ ghost points at interval ends. This gives $(r_i+K)^d$ basis functions. Figure~\ref{fig:B-spline} illustrates a 1d B-spline basis. The function ansatz is
\begin{equation}
    f(\vec x) = \sum_{n_1=1}^{r_1+K} \dots \sum_{n_d=1}^{r_d+K} \mathcal{C}_{n_1, \dots, n_d}\, \rho_1(x^{(1)})\dots\rho_d(x^{(d)}),
\end{equation}
where $\mathcal{C}$ is the coefficient tensor, which can be sampled from uniform or normal distributions. If the coefficients lie in $[a,b]$, then $\|f\|_{L^\infty(\Omega)} \le \max(|a|,|b|)$, ensuring the valid input range of the model. For example, this can be used to guarantee that PDE coefficients exceed a given threshold.

For simple convex domains, the spline domain can inscribe the region of interest. More complex geometries may require a specialized ansatz, such as a spline basis defined on a finite element mesh.

\subsection{Input scaling}
Proper input scaling is crucial for neural network training. Simply drawing the coefficient tensor from a distribution may not produce function values with a balanced distribution across sensor points. With B-splines, the core tensor can be sampled within a defined range. A more effective approach is to normalize the input to have zero mean and unit variance. In higher dimensions, computing variance across sensor points can be memory-intensive if done naively. This can be addressed efficiently using a single-pass algorithm such as Welford's method \cite{welford1962note}, where functions or batches of functions are processed sequentially, updating the variance estimate incrementally.

\subsection{Penalty framework}
If we consider for example a time independent PDE of the form
\begin{equation}
    \begin{aligned}
        \mathcal{L}(u)(\vec x) &= f(\vec x) \quad \text{for}\, \vec x \in \Omega \\
        u(\vec x) &= 0 \quad \text{for}\, \vec x \in \partial\Omega,
    \end{aligned}
\end{equation}
with differential operator $\mathcal{L}$ and homogeneous Dirichlet boundary conditions,
we use our LNO model $\mathcal{M}$ with parameters $\Theta$ as an ansatz for the solution $u$ such that after training $u(\vec Y) \approx \mathcal{M}[\vec Y, \vec X, f, \Theta]$ . Now we can sample $i=1,\dots,p$ source term functions $f_i$ from our functional space and evaluate them at sensor locations $\vec X_i$ to further make predictions on some target points $\vec Y_i \subset \Omega$ as well as target points on the boundary $\vec{\bar{Y}}_i \subset \partial\Omega$. Here we term $p$ the batch size and note that for each sample, sensor and target points could be different. However, we find that using the same points across all samples also works reasonably well for our purposes.
We can use a physics informed loss
\begin{equation}
    J_{\text{PDE}}(\Theta) = \frac{1}{p}\sum_{i=1}^p \frac{1}{|\vec Y_i|}\sum_{y \in \vec Y_i} \big(\mathcal{L}(\mathcal{M}(\vec X_i, f_i, \Theta))(y) - f_i(\vec x)\big)^2.
\end{equation}
Boundary conditions can be accounted for with a penalty term
\begin{equation}
    J_{\text{B}}(\Theta) = \frac{1}{p}\sum_{i=1}^p \frac{1}{|\vec{\bar{Y}}_i|}\sum_{y \in \vec{\bar{Y}}_i} \mathcal{M}(\vec X_i, f_i, \Theta)(y)^2.
\end{equation}
The physics informed loss is then given by $J_{\text{PI}} = J_{\text{PDE}} + \lambda J_{\text{B}}$, where $\lambda$ is a penalty parameter. This empirical risk can be minimized with respect to $\Theta$ by means of stochastic gradient descent. Since we can sample unlimited data instances, we don't speak from epochs, as is usually the case in machine learning, but simply from steps. It is not necessary to resample $\vec X_i, \vec{\bar X}_i$ and $\vec Y_i$ for each step, but it is sufficient to resample if after a certain number of steps has passed. 
It is a good strategy to also resample the sensor points, since that allows our model to perform well even on unseen sensor points. In our setting, only the functions $f_i$ are newly sampled for each step and the penalty parameter is gradually increased during training once a predefined number of steps is reached. This penalty framework might be heavy, but since we are learning the solution operator to a whole class of PDEs and not just a single PDE as is the case for PINNs, we can afford to spend larger amounts of resources on training. This section should serve as a general guideline to train such operators, and it should be possible to adapt this scheme to many different problem settings. 

\section{Results}
To demonstrate our LNO model, we show its performance on some toy problems. First, we consider simple function fitting with our LNO model. Afterwards, we use the physics-informed penalty framework to train an LNO model for elliptic PDEs. Elliptic PDEs form a cornerstone of mathematical modeling, capturing steady-state processes at equilibrium. Classical examples include Laplace’s and Poisson’s equations, which describe electrostatics, stationary heat conduction, gravitation, or static elasticity. Therefore, we use variants of elliptic PDEs to further test our LNO model.

Our computational domain is $\Omega=[-0.5, 0.5]^2$. Sensor points $\vec X$ and target points $\vec Y$ are drawn from a Sobol sequence \cite{sobol1976uniformly} with $N=2^{10}$ sensor points and $M=2^9$ target points during training. Boundary target points $\bar{\vec Y}$ are drawn uniformly with $\bar M=2^8$ points. These points are updated every $\num{500}$ steps. For all experiments, the batch size is $p=30$, with source and target points identical across instances. 

Function sampling uses a B-spline basis with $K=3$ and $r_1 = r_2 = 10$, yielding a multidimensional rank of $\num{13}$. The coefficient tensor is drawn from a standard normal distribution. We denote an MLP with hidden layers $[a_1,\dots,a_P]$, where $a_P$ is the output dimension. The input dimension is omitted. Each layer includes a bias term, and a \texttt{tanh} activation is used for all layers except the output layers of $\mathcal{N}_T$, $\psi_t$, and $\phi_t$ for $t=1,\dots,T$, which have linear output activations. All layers are initialized using Xavier normal initialization.

\subsection{Function fitting}

As a first toy example, we consider simple function fitting with our LNO model. The function values are given at the sensor points, and the model aims to create a continuous representation by encoding the input into a latent space. The model uses $T=4$ LNO layers, an embedding dimension of $S=50$, and a low-rank dimension $R=100$. The embedding network $\mathcal{N}_0$ has layers $[S, S, S]$, while $\mathcal{N}_t$ for $t=1,\dots, T-1$ has layers $[S, S]$, and the decoder $\mathcal{N}_T$ has layers $[S, 10, 1]$. A symmetric kernel ($\psi_t=\phi_t$ for $t=1,\dots,T$) is used, and each $\phi_t$ has layers $[50, 50, R]$.

\begin{figure}
    \centering
    \includegraphics[width=1.0\linewidth]{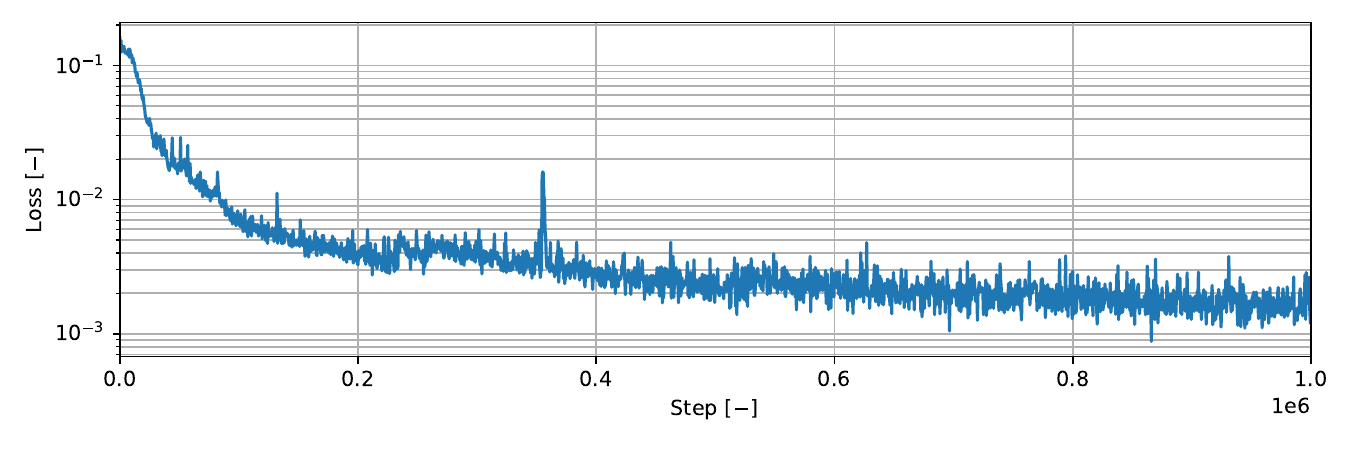}
    \caption{Learning curve of the LNO model for function fitting.}
    \label{fig:lc-function-fitting}
\end{figure}

The loss function for a mini-batch $\{f_1,\dots,f_p\}$ is
\begin{equation}\label{eq:loss-fn-fitting}
    J_{\text{MSE}}(\Theta) = \frac{1}{p}\sum_{i=1}^p \frac{1}{M}\sum_{y \in \vec{Y}} \big(\mathcal{M}(\vec X, f_i, \Theta)(y) - f_i(y)\big)^2.
\end{equation}
The model is trained for $\num{1000000}$ steps with a learning rate of $\num{5e-4}$. Figure~\ref{fig:lc-function-fitting} shows the loss function during training.

\begin{figure}
    \centering
    \includegraphics[width=1\linewidth]{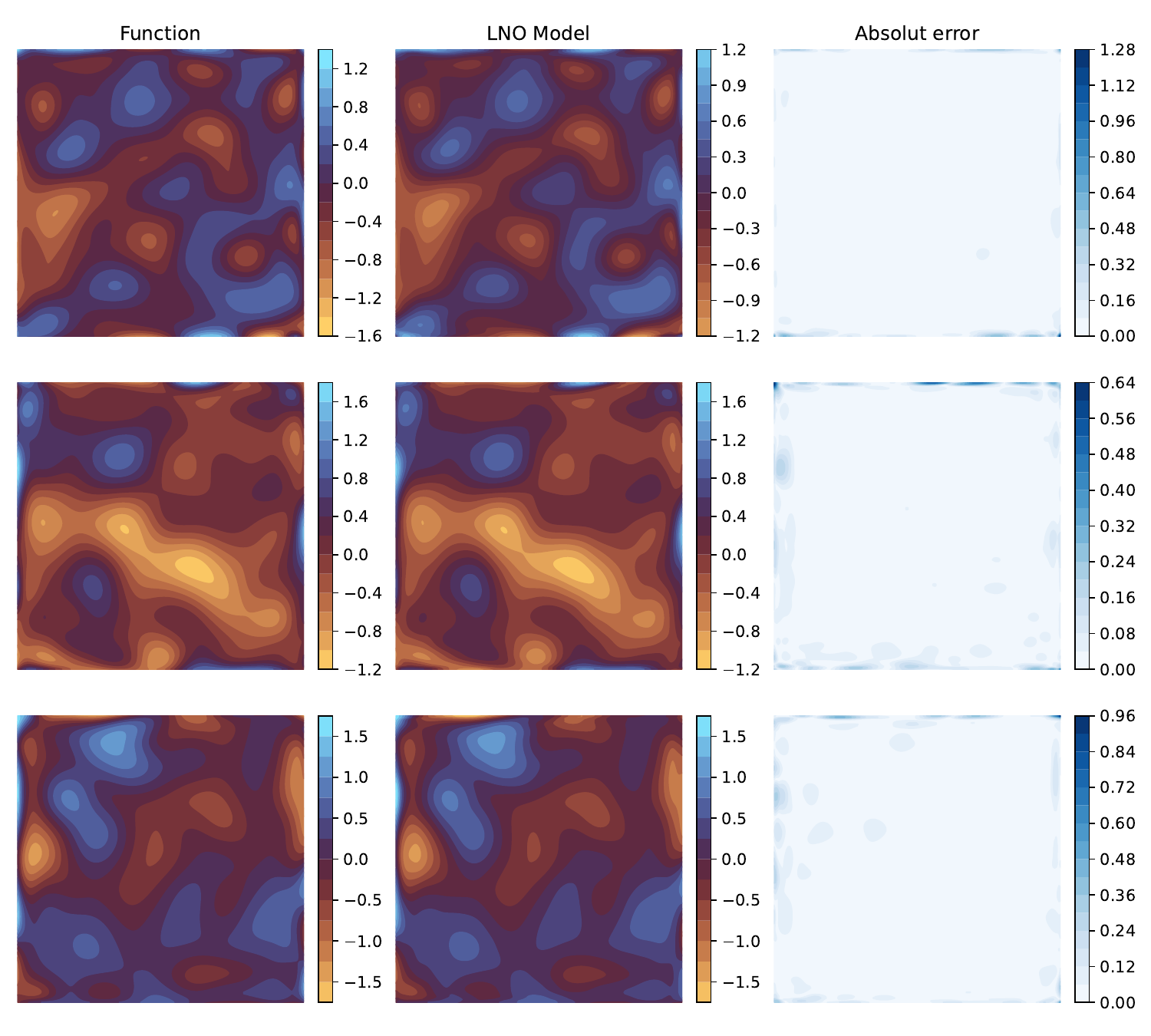}
    \caption{Function fitting with an LNO model using $N=2^{12}$ sensor points.}
    \label{fig:function-fitting}
\end{figure}

Figure~\ref{fig:function-fitting} shows the model predictions for some unseen functions. The model performs well overall but shows errors near the domain boundaries.

\begin{table}[h]
    \centering
    \begin{tabular}{c| c | c | c | c | c | c | c | c }
        N & $64$ & $128$ &$256$ &$512$ & $1024$ & $2048$ & $4096$ & $8192$ \\
        \hline
        RMSE & $\num{0.797}$ & $\num{0.435}$ & $\num{0.206}$ & $\num{0.123}$ & $\num{0.089}$ & $\num{0.052}$ & $\num{0.050}$ & $\num{0.047}$ \\
        GPU pred. time $[\si{\milli\second}]$ & $\num{0.262}$ & $\num{0.256}$ & $\num{0.212}$ & $\num{0.221}$ & $\num{0.250}$ & $\num{0.227}$ & $\num{0.225}$ & $\num{0.208}$ \\
        CPU pred. time $[\si{\milli\second}]$ & $\num{1.872}$ & $\num{2.065}$ & $\num{2.146}$ & $\num{2.065}$ & $\num{2.381}$ & $\num{2.182}$ & $\num{2.939}$ & $\num{4.676}$ \\
    \end{tabular}
    \caption{Mean prediction error and mean prediction time for GPU and CPU for the trained function fitting LNO model over $p=\num{1000}$ samples for a varying number of sensor points $N$ and $M=\num{10000}$ target points. }
    \label{tab:fn-fitting-table}
\end{table}
Table~\ref{tab:fn-fitting-table} shows the root mean squared error (RMSE) and prediction time on GPU and CPU for different numbers of sensor points. The RMSE is computed as the square root of \eqref{eq:loss-fn-fitting} using the trained model parameters. A small number of sensor points leads to poor performance, while the prediction time on GPU remains essentially constant due to parallelism.

\subsection{Poisson equation}
As a second test, we approximate the solution operator of the Poisson equation with homogeneous Dirichlet boundary conditions on $\Omega = [-0.5, 0.5]^2$:
\begin{equation}
\begin{aligned}
    -\Delta u &= f, \quad x \in \Omega \\
    u &= 0, \quad x \in \partial\Omega.
\end{aligned}
\end{equation}
The Poisson equation is a fundamental elliptic PDE, modeling how a potential field $u$ arises from a source distribution $f$.

For our LNO model, we use embedding size $S=50$ and low-rank $R=200$, with $T=7$ layers and a symmetric kernel. Each $\phi_t$ for $t=1,\dots,T$ is a network with layers $[50, 50, 50, R]$. The embedding $\mathcal{N}_0$ has layers $[S, S, S]$, $\mathcal{N}_t$ for $t=1,\dots,T-1$ has layers $[S, S]$, and the decoder $\mathcal{N}_T$ has layers $[S, 10, 1]$, resulting in $\num{147623}$ trainable parameters.

\begin{figure}
    \centering
    \includegraphics[width=1\linewidth]{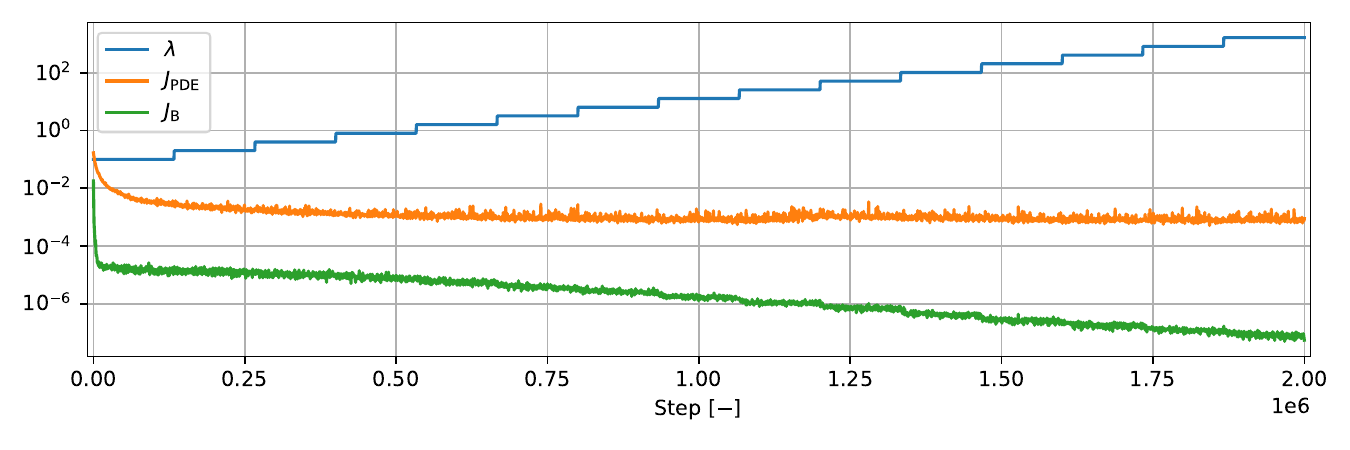}
    \caption{Penalty framework training for the Poisson equation. At the end of training $J_{\text{PDE}}$ is around $\num{9.5e-4}$ and $J_{\text{B}}$ around $\num{6.7e-8}$.}
    \label{fig:poisson-lc}
\end{figure}

The source is given by a B-spline ansatz, and the physics-informed loss $J_{\text{PI}} = J_{\text{PDE}} + \lambda J_{\text{B}}$ is minimized using a penalty framework. The penalty parameter $\lambda$ starts at $\lambda_0=0.1$ and is gradually doubled during training, as seen in the learning curve in Figure~\ref{fig:poisson-lc}. The model was trained with the Muon optimizer \cite{jordan2024muon} for $\num{2000000}$ steps with learning rate $\num{5e-4}$. Earlier attempts with Adam did not yield satisfactory results. Muon, which incorporates curvature information, outperformed Adam in all experiments. The model is still not fully optimized, and better optimization schemes, possibly stochastic second-order methods, could further improve performance.

\begin{figure}
    \centering
    \includegraphics[width=1\linewidth]{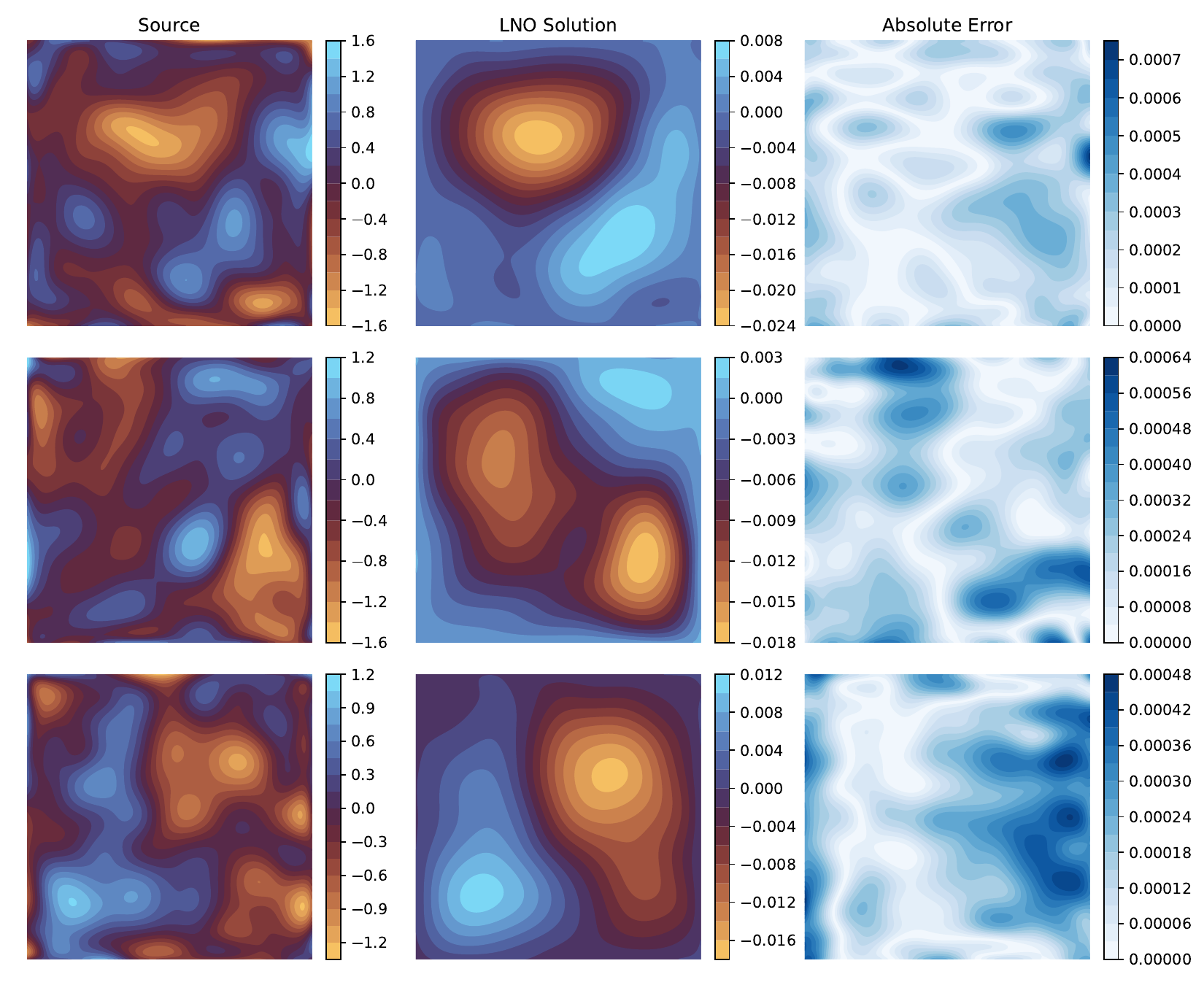}
    \caption{Prediction of the LNO model for some sampled source functions from the training space, together with the absolute error with respect to a finite difference solution on a grid with $250\times 250$ cells. The model used $N=1024$ sensor points for the prediction.}
    \label{fig:poisson}
\end{figure}

Figure~\ref{fig:poisson} shows predictions of the LNO model for unseen source functions with $N=1024$ sensor points.

\begin{figure}
    \centering
    \includegraphics[width=0.9\linewidth]{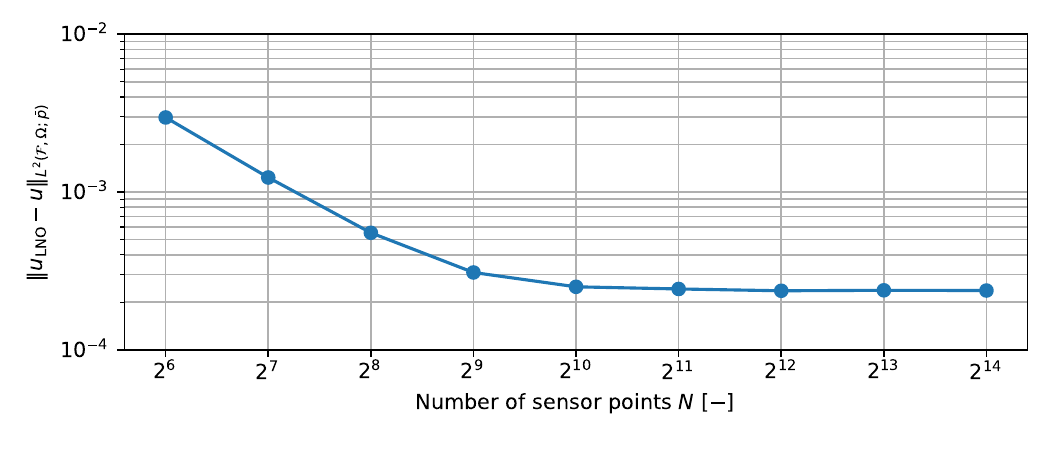}
    \caption{Prediction error of the LNO model for the Poisson equations a for varying number of sensor points.}
    \label{fig:poisson-error}
\end{figure}

One important question is the model performance with varying numbers of sensor points. To investigate this, we increase the number of sensor points and compute the solution for $\bar p = 500$ different source terms. Figure~\ref{fig:poisson-error} shows the empirical $L^2$ error over these samples with respect to the finite difference solution:
\begin{equation}
    \|u_{\text{LNO}} - u\|_{{L^2}(\mathcal{F}, \Omega ; p)} = \frac{1}{p}\sum_{i=1}^p \|u_{\text{LNO}} - u\|_{{L^2}(\Omega)}.
\end{equation}
With too few sensor points, the accuracy decreases. However, using more sensor points than in training does not further improve performance. The average relative $L_2$ error is approximately $5\%$ for a sufficient number of sensor points.

\begin{figure}
    \centering
    \includegraphics[width=1\linewidth]{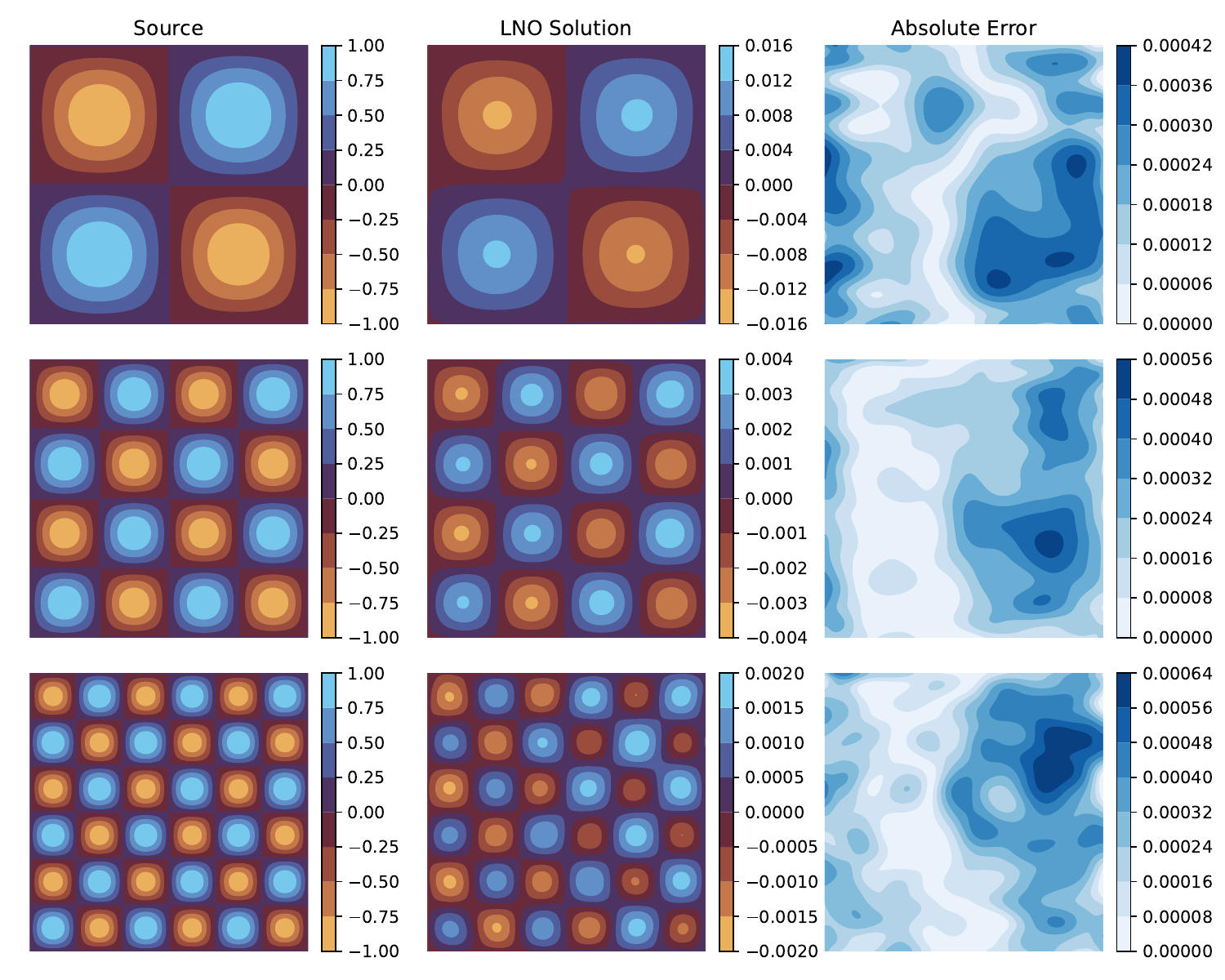}
    \caption{Prediction of the Poisson LNO model for the source $f(x_1, x_2) = \sin(k \pi x_1) \sin(k \pi x_2)$ with $k=2$ (top), $k=4$ (middle) and $k=6$ (bottom) and $N=2^{10}$ sensor points.}
    \label{fig:poisson-sin}
\end{figure}

To show, that the model indeed learns a useful kernel, we predict the solution to a more structured problem given by the source
\begin{equation}\label{eq:sin-source}
    f(x_1, x_2) = \sin(k \pi x_1) \sin(k \pi x_2),
\end{equation}
for $k=2,4,6$, with the solution given by
\begin{equation}
    u(x_1, x_2) = \frac{1}{2 k^2 \pi^2}\sin(k \pi x_1) \sin(k \pi x_2).
\end{equation}
Figure~\ref{fig:poisson-sin} shows the model prediction for this source.

\subsection{Screened Poisson}
The screened Poisson equation is a modification of the classical Poisson equation and is given by
\[
-\Delta u(x) + s u(x) = f(x), \quad x \in \Omega,
\]
where the parameter $s > 0$ represents a damping or screening effect. This equation appears in various applications, including electrostatics with Debye–Hückel screening, quantum mechanics (Yukawa potential), and image or geometry processing for reconstruction and smoothing tasks. Unlike the classical Poisson equation, the additional term $s u$ ensures that the influence of sources $f$ decays locally rather than propagating globally.

\begin{figure}
    \centering
    \includegraphics[width=1\linewidth]{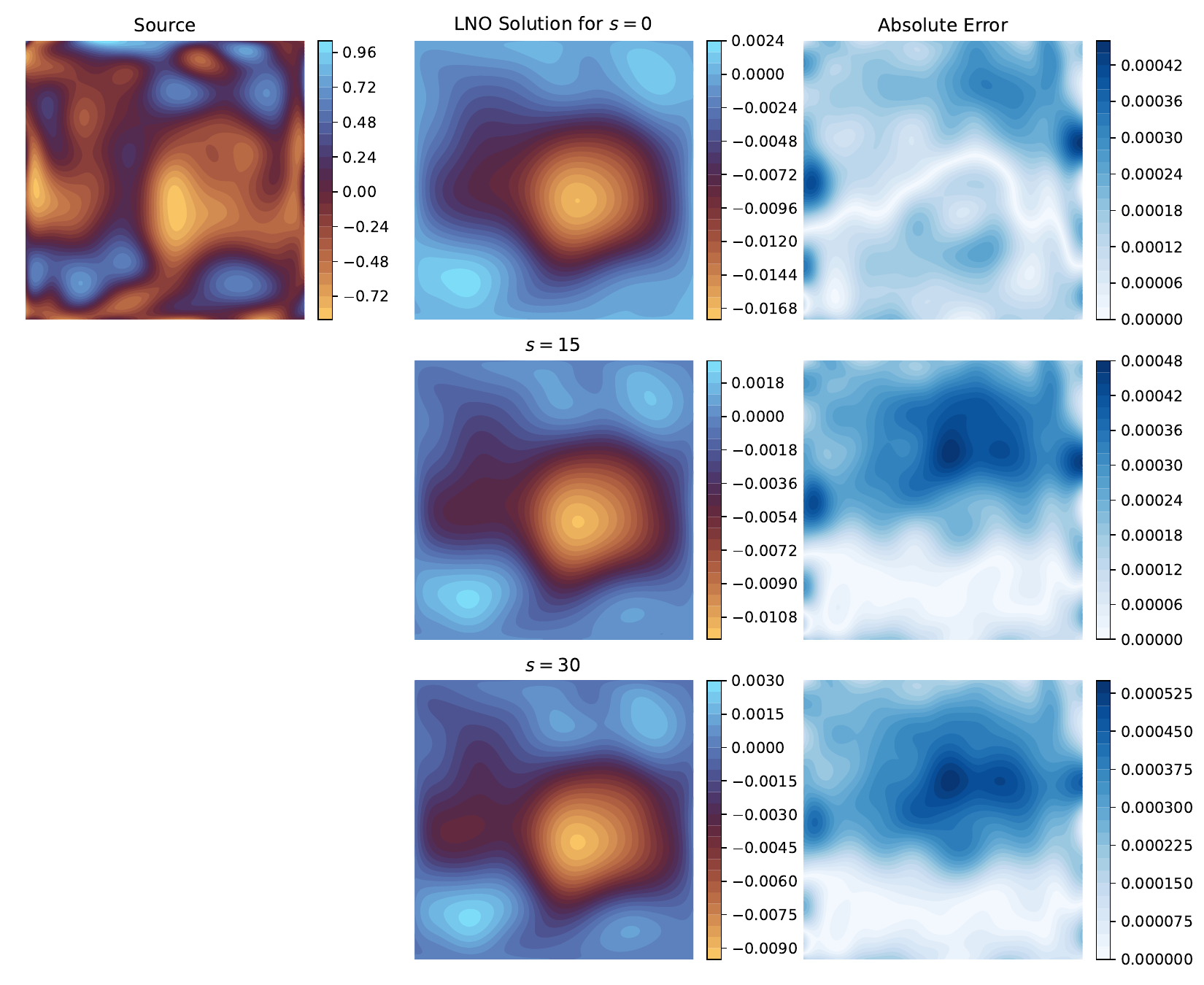}
    \caption{Prediction of the screened Poisson LNO model for a single source and different values of the screening coefficient. The absolute error was computed with respect to a finite difference solution on a grid with $250\times 250$ cells}
    \label{fig:screened-poisson}
\end{figure}

In our experiments, we approximate the solution operator for the screening coefficient $s\in[0,30]$ using the PILNO model. The model is essentially the same as for the Poisson problem, with the difference that the symmetric low-rank kernel functions also take $s$ as input, i.e.,
\[
k_t(x, y, s) = \phi_t(x, s / 30)^T \phi_t(y, s / 30), \quad t=1,\dots,T,
\]
and we use $T=7$ layers. During training, $s$ is drawn uniformly from its domain together with the source function $f$, giving a training batch $\{(s_1, f_1), \dots, (s_p, f_p)\}$. The training curves for the penalty framework are nearly identical to those in Figure~\ref{fig:poisson-lc}. The final loss is $J_{\text{PDE}}=\num{6.5e-4}$ and $J_{\text{B}}=\num{4.2e-8}$. Although the loss is slightly better than for the Poisson model, the average relative $L^2$ error is around 7.5\% for a sufficient number of sensor points. Figure~\ref{fig:screened-poisson} shows the prediction for a specific source with three different values $s \in \{0, 15, 30\}$.

This LNO model learns a whole family of solution operators based on the conditional parameter $s$.

\subsection{Darcy flow}
The Darcy flow equation describes fluid movement through porous media and is formulated as
\begin{equation}
\begin{aligned}
    -\nabla \cdot \left(c(x) \nabla u(x)\right) &= f \quad \text{for } x \in \Omega,
\end{aligned}
\end{equation}
where $u$ denotes the fluid pressure, $c$ is the (parameter-dependent) diffusion coefficient of the medium, and $f$ represents sources or sinks such as injection or extraction wells. Darcy flow is important in hydrogeology, petroleum engineering, and environmental sciences, with applications ranging from groundwater management and oil recovery to contaminant transport and carbon dioxide sequestration. Its relevance stems from capturing how heterogeneous porous structures influence large-scale flow patterns. Homogeneous Dirichlet boundary conditions are used.

In this example, we train the model for a high-dimensional parameter space. Using a $\bar K$-dimensional basis expansion for the diffusion coefficient:
\begin{equation}
    c(x) = \sum_{i=1}^{\bar K}\beta_i \rho_i(x),
\end{equation}
this represents a training task over the $\bar K$-dimensional parameter space plus the problem dimension. Conventional numerical methods require assembling the linear system for each choice of $c$, which is infeasible for high-dimensional spaces. Neural operators, however, can approximate the solution manifold efficiently.

The model adjusts to this high-dimensional setting by including the diffusion coefficient pointwise in the kernel:
\[
k_t(x, y, c) = \psi_t(x, c(x))^T\phi_t(y, c(y)), \quad t=1,\dots,T,
\]
and since the solution operator is not necessarily symmetric, we use $\psi_t \neq \phi_t$. For unsupervised training, B-splines are used to sample $c$, ensuring $\beta_i > 0$ for $i=1,\dots,\bar K$ so that $c(x) > 0$ in $\Omega$. Training batches are $\{(f_1, c_1), \dots, (f_p, c_p)\}$. Training was the same as for the Poisson problem, with final losses $J_{\text{PDE}}=\num{1.23e-4}$ and $J_{\text{B}} \approx \num{1e-7}$.

In particular, we sample the core tensor such that $\beta_i \in (0.2, 1)$ for $i=1,\dots,\bar K$ and use a B-spline with $r_1=r_2=5$ and order $3$, resulting in $\bar K=64$. This B-spline is used for both the source function and the diffusion coefficient. We use the same architecture as for the screened Poisson problem (except that $\psi_t$ and $\phi_t$ do not share parameters) and note that the function space for the Darcy flow is slightly less complex. Small function values of $c$ are avoided, as they lead to sharp solution gradients, which are challenging for both PINNs and neural operators. For illustration, this simpler setting suffices.
\begin{figure}
    \centering
    \includegraphics[width=1\linewidth]{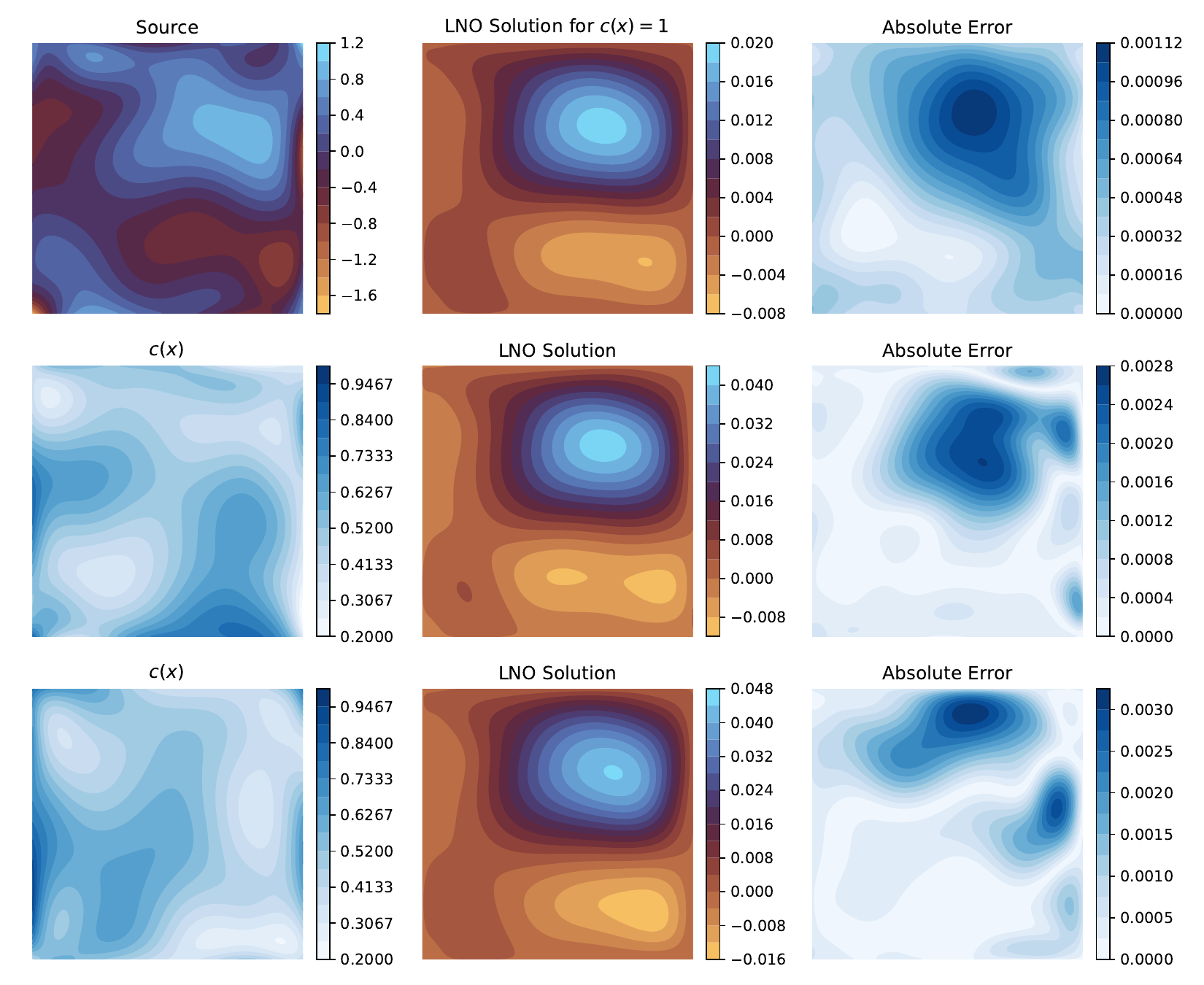}
    \caption{LNO predictions for the Darcy flow. The top row shows the source $f$, the prediction of the model for constant diffusion coefficient and the absolute error with respect to a finite difference solution. Middle and bottom row show the prediction for the same source but different coefficient functions.}
    \label{fig:diffusion}
\end{figure}

Figure~\ref{fig:diffusion} shows predictions of the LNO model for a sampled source function and different diffusion coefficients. Overall, the model achieves a mean relative $L_2$ error of about $\num{14}\%$, computed over 100 samples of sources and coefficients.

\begin{figure}
    \centering
    \includegraphics[width=1\linewidth]{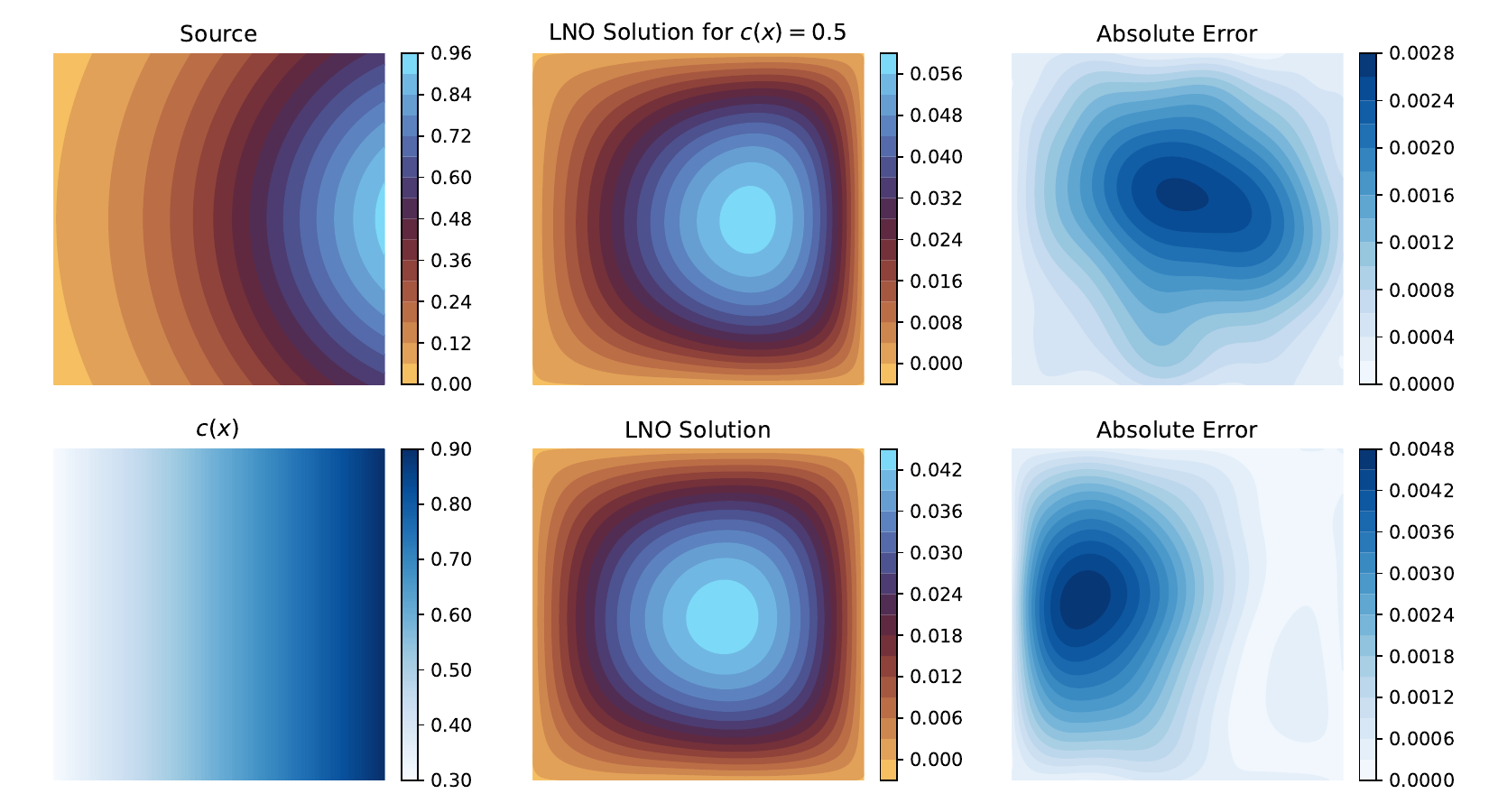}
    \caption{LNO prediction for the source $f(x)=\exp(-2 \|x - [0.7, 0]^T\|^2)$ and constant diffusion coefficient (top row) and for the coefficient function $c(x) = 0.6x^{(1)} + 0.6$ (bottom row).}
    \label{fig:diffusion2}
\end{figure}
Figure~\ref{fig:diffusion2} shows predictions for a more structured example, demonstrating that random B-spline sampling efficiently covers the solution space. Specifically, we use $f(x)=\exp(-2 \|x - [0.7, 0]^T\|^2)$ and $c(x) = 0.6 x^{(1)} + 0.6$.
While the model is not yet accurate enough for practical applications, optimization improvements could enhance performance. For example, energy natural gradient descent \cite{muller2023achieving} improves PINN solutions, but direct application to stochastic neural operator training is nontrivial. Momentum-based extensions \cite{guzman2025improving} or Kronecker-factored approximate curvature methods \cite{dangel2024kronecker} may therefore be promising.

\section{Conclusion}
We introduced the Physics-Informed Low-Rank Neural Operator (PILNO), a neural operator
architecture for learning solution operators of partial differential equations (PDEs) on point cloud data. PILNO combines low-rank kernel approximations with an encoder--decoder
architecture, enabling efficient approximation of convolution-like operations in
high-dimensional parameter spaces while remaining independent of specific discretizations.
The model was validated on a range of problems, including function fitting, the Poisson
equation, the screened Poisson equation with variable coefficients, and parameterized Darcy flow.

Our results demonstrate that PILNO can deliver continuous, one-shot predictions with a
limited number of sensor points. The low-rank structure ensures computational efficiency
during both encoding and decoding, and the physics-informed penalty framework enforces
PDE constraints and boundary conditions, even without labeled data. Although classical
mesh-based solvers currently provide higher absolute accuracy, PILNO offers significant
advantages in efficiency and scalability, particularly for high-dimensional parameter spaces. 
In fact, PILNO trades some accuracy for efficiency.
Improving kernel parameterization and training stability remains a key challenge.

Looking ahead, several avenues for further research are promising. These include the
development of advanced and accurate optimization strategies for training neural operators, the extension to more complex geometries and boundary conditions, and the incorporation of problem-specific features into kernel parameterizations. Beyond static elliptic problems,
PILNO could be extended to time-dependent PDEs, multi-physics systems, and adaptive sampling strategies to enhance accuracy and generalization. Together, these directions point to PILNO as a versatile and scalable surrogate modeling framework that bridges classical numerical methods and modern data-driven approaches.

\section*{Acknowledgements}
\noindent Financial support by the Austrian Science Fund (FWF) via project ”Data-driven Reduced Order Approaches for Micromagnetism (Data-ROAM)” (Grant-DOI: 10.55776/PAT7615923), project ”Design of Nanocomposite Magnets by
Machine Learning (DeNaMML)” (Grant-DOI: 10.55776/P35413) is gratefully acknowledged. The authors acknowledge the University of Vienna research platform MMM Mathematics - Magnetism - Materials. The computations were partly achieved by using the Vienna Scientific Cluster (VSC) via the funded projects No. 71140, 71952 and 72862.
This research was funded in whole or in part by the Austrian Science Fund (FWF) [10.55776/PAT7615923, 10.55776/P35413]. For the purpose of Open Access, the authors have applied a CC BY public copyright license to any Author Accepted Manuscript (AAM) version arising from this submission. 


\begin{thebibliography}{10}

\bibitem{raissi2019physics}
M.~Raissi, P.~Perdikaris, and G.~E. Karniadakis, ``Physics-{I}nformed {N}eural {N}etworks: A deep learning framework for solving forward and inverse problems involving nonlinear partial differential equations,'' {\em Journal of Computational physics}, vol.~378, pp.~686--707, 2019.

\bibitem{kovacs2022conditional}
A.~Kovacs, L.~Exl, A.~Kornell, J.~Fischbacher, M.~Hovorka, M.~Gusenbauer, L.~Breth, H.~Oezelt, M.~Yano, N.~Sakuma, {\em et~al.}, ``Conditional physics informed neural networks,'' {\em Communications in Nonlinear Science and Numerical Simulation}, vol.~104, p.~106041, 2022.

\bibitem{lu2021learning}
L.~Lu, P.~Jin, G.~Pang, Z.~Zhang, and G.~E. Karniadakis, ``Learning nonlinear operators via deeponet based on the universal approximation theorem of operators,'' {\em Nature machine intelligence}, vol.~3, no.~3, pp.~218--229, 2021.

\bibitem{li2024physics}
Z.~Li, H.~Zheng, N.~Kovachki, D.~Jin, H.~Chen, B.~Liu, K.~Azizzadenesheli, and A.~Anandkumar, ``Physics-informed neural operator for learning partial differential equations,'' {\em ACM / IMS J. Data Sci.}, vol.~1, May 2024.

\bibitem{wang2021learning}
S.~Wang, H.~Wang, and P.~Perdikaris, ``Learning the solution operator of parametric partial differential equations with physics-informed {D}eep{O}{N}ets,'' {\em Science advances}, vol.~7, no.~40, p.~eabi8605, 2021.

\bibitem{li2020neural}
Z.~Li, N.~Kovachki, K.~Azizzadenesheli, B.~Liu, K.~Bhattacharya, A.~Stuart, and A.~Anandkumar, ``Neural operator: Graph kernel network for partial differential equations,'' {\em arXiv preprint arXiv:2003.03485}, 2020.

\bibitem{li2020fourier}
Z.~Li, N.~Kovachki, K.~Azizzadenesheli, B.~Liu, K.~Bhattacharya, A.~Stuart, and A.~Anandkumar, ``Fourier neural operator for parametric partial differential equations,'' {\em arXiv preprint arXiv:2010.08895}, 2020.

\bibitem{li2023fourier}
Z.~Li, D.~Z. Huang, B.~Liu, and A.~Anandkumar, ``Fourier neural operator with learned deformations for pdes on general geometries,'' {\em Journal of Machine Learning Research}, vol.~24, no.~388, pp.~1--26, 2023.

\bibitem{kossaifi2023multi}
J.~Kossaifi, N.~Kovachki, K.~Azizzadenesheli, and A.~Anandkumar, ``Multi-grid tensorized fourier neural operator for high-resolution pdes,'' {\em arXiv preprint arXiv:2310.00120}, 2023.

\bibitem{kovachki2023neural}
N.~Kovachki, Z.~Li, B.~Liu, K.~Azizzadenesheli, K.~Bhattacharya, A.~Stuart, and A.~Anandkumar, ``Neural operator: Learning maps between function spaces with applications to pdes,'' {\em Journal of Machine Learning Research}, vol.~24, no.~89, pp.~1--97, 2023.

\bibitem{zeng2025point}
C.~Zeng, Y.~Zhang, J.~Zhou, Y.~Wang, Z.~Wang, Y.~Liu, L.~Wu, and D.~Z. Huang, ``Point cloud neural operator for parametric pdes on complex and variable geometries,'' {\em Computer Methods in Applied Mechanics and Engineering}, vol.~443, p.~118022, 2025.

\bibitem{sukumar2022exact}
N.~Sukumar and A.~Srivastava, ``Exact imposition of boundary conditions with distance functions in physics-informed deep neural networks,'' {\em Computer Methods in Applied Mechanics and Engineering}, vol.~389, p.~114333, 2022.

\bibitem{schaffer2024constraint}
S.~Schaffer and L.~Exl, ``Constraint free physics-informed machine learning for micromagnetic energy minimization,'' {\em Computer Physics Communications}, vol.~300, p.~109202, 2024.

\bibitem{exl2025higher}
L.~Exl and S.~Schaffer, ``Higher order stray field computation on tensor product domains,'' {\em arXiv preprint arXiv:2505.19180}, 2025.

\bibitem{welford1962note}
B.~P. Welford, ``Note on a method for calculating corrected sums of squares and products,'' {\em Technometrics}, vol.~4, no.~3, pp.~419--420, 1962.

\bibitem{sobol1976uniformly}
I.~M. Sobol, ``Uniformly distributed sequences with an additional uniform property,'' {\em USSR Computational Mathematics and Mathematical Physics}, vol.~16, no.~5, pp.~236--242, 1976.

\bibitem{jordan2024muon}
K.~Jordan, Y.~Jin, V.~Boza, J.~You, F.~Cesista, L.~Newhouse, and J.~Bernstein, ``Muon: An optimizer for hidden layers in neural networks,'' 2024.

\bibitem{muller2023achieving}
J.~M{\"u}ller and M.~Zeinhofer, ``Achieving high accuracy with {PINN}s via energy natural gradient descent,'' in {\em International Conference on Machine Learning}, pp.~25471--25485, PMLR, 2023.

\bibitem{guzman2025improving}
A.~Guzm{\'a}n-Cordero, F.~Dangel, G.~Goldshlager, and M.~Zeinhofer, ``Improving energy natural gradient descent through woodbury, momentum, and randomization,'' {\em arXiv preprint arXiv:2505.12149}, 2025.

\bibitem{dangel2024kronecker}
F.~Dangel, J.~M{\"u}ller, and M.~Zeinhofer, ``Kronecker-factored approximate curvature for physics-informed neural networks,'' {\em Advances in Neural Information Processing Systems}, vol.~37, pp.~34582--34636, 2024.

\end{thebibliography}

\end{document}